\documentclass[12pt,oneside]{amsart}
\usepackage{amsmath, amsfonts,amsthm,times,graphics}
\input amssym.def
\input amssym

\begin{document}

\newtheorem{theorem}{Theorem}[section]
\newtheorem{lemma}[theorem]{Lemma}
\newtheorem{conjecture}[theorem]{Conjecture}
\newtheorem{question}[theorem]{Question}
\newtheorem{problem}[theorem]{Problem}
\newtheorem{corollary}[theorem]{Corollary}
\newtheorem{proposition}[theorem]{Proposition}

\newtheorem{definition}[theorem]{Definition}
\newtheorem{xca}[theorem]{Exercise}
\theoremstyle{remark}
\newtheorem{remark}[theorem]{Remark}
\newtheorem{example}[theorem]{Example}

\title[Goldbach's like conjectures 
arising from arithmetic progressions\ldots] 
{Goldbach's like conjectures arising from arithmetic progressions
 whose first two terms are primes}

\author{Romeo Me\v strovi\' c}
\address{Maritime Faculty Kotor, University of Montenegro, Dobrota 36,
85330 Kotor, Montenegro} \email{romeo@ucg.ac.me}

\makeatother

{\renewcommand{\thefootnote}{}
\setcounter{footnote}{0}}

\maketitle

 \begin{abstract}
For two  odd primes  $p$ and $q$ such that $p<q$, let  
$A(p,q):=(a_k)_{k=1}^{\infty}$  be the arithmetic progression 
whose $k$th term is given by $a_k=(k-1)(q-p)+p$ (i.e., with  $a_1=p$ and $a_2=q$).
Here we conjecture that for every positive integer $a>1$
there exist a positive integer $n$ 
and  two odd primes $p$ and $q$ such that 
$a$ can be expressed as a sum of the first $2n$ terms of 
the arithmetic progression $A(p,q)$.
 Notice that in the case of even $a$,
this conjecture immediately follows from Goldbach's conjecture.
We also propose the analogous conjecture for odd  positive integers $a>1$
as well as some related Goldbach's like conjectures arising 
from the   previously mentioned arithmetic progressions.
 \end{abstract}
\section{Conjectures
  on  arithmetic progressions whose first two terms are primes}

Let $p$ and $q$ be two  primes such that $p<q$ and let  
$A(p,q):=(a_k)_{k=1}^{\infty}$ be the arithmetic progression 
whose $k$th term is given by 
  $$
a_k=(k-1)(q-p)+p, \quad k=1,2,\ldots
  $$
 In other words, $A(p,q)$ is 
an arithmetic progression whose first two terms are $p$ and $q$
(i.e., $a_1=p$ and $a_2=q$).
The sum $S_n(p,q)=S_n$ of the first $n$ terms of 
the progression $A(p,q)$ is equal to
  $$
S_{n}(p,q)=\frac{n}{2}((n-1)q-(n-3)p).\leqno(1) 
  $$
From (1) we have that for all 
$n=1,2,\ldots$ and $m=0,1,2,\ldots$ the sum $S_{n,m}(p,q):=\sum_{i=m+1}^{n+m}a_i$
of some $n$ consecutive terms of progression $A(p,q)$ is equal to
 $$
S_{n,m}(p,q):=S_{n+m}(p,q)-S_{m}(p,q)=
\frac{n}{2}((n+2m-1)q-(n+2m-3)p).\leqno(2) 
  $$

We start with following example.

  \begin{example}[An extension of a Sylvester's result]
Here we examine positive integers $a$ which can be written 
as a sum $S_{n,m}(2,3)$  (given by (2) with $p=2$ and $q=3$) for some 
$n\ge 2$ and $m\ge 1$. The sum of $k$th term and $(k+1)$th term 
of the progression $A(2,3)=(k+1)_{k=1}^{\infty}$ is equal to $2k+3$.
Therefore, every odd integer greater than 3 is a sum of 
some two  consecutive terms of $A(2,3)$. Furthermore, by (2) we have
  $$
S_{n,m}(2,3)=\sum_{i=m+1}^{n+m}(i+1)=\frac{n}{2}(n+2m+3).\leqno(3) 
  $$
If $a$ is an even positive integer which is not a power of 2, 
then $a=(2d+1)2^u$ for some positive integers $d\ge 1$ and $u\ge 1$. If 
$1\le d\le 2^u-2$ for such a $a$, we have  
 $S_{2d+1,2^u-d-2}=(2d+1)2^u=a$.
(If $d=0$ then $n=1$ and $S_{1,m}(2,3)=m+2$ is in fact the $(m+1)$th term 
of $A(2,3)$).  
Similarly, if $d\ge 2^u+1$, then $S_{2^{u+1},d-2^u-1}=(2d+1)2^u=a$.
This shows that each even positive integer $a=(2d+1)2^u$
with  $1\le d\le 2^u-2$ or $d\ge 2^u+1$ can be expressed as a sum 
of at least two consecutive  terms of the arithmetic progression $A(2,3)$.

It remains to consider the cases when $a$ is of the form  
$2^u$, $(2^{u+1}-1)2^u$ or $(2^{u+1}+1)2^u$ with some positive integer $u$. 
If $a=2^u$, then by (3) the equality $S_{n,m}(2,3)=a$ is equivalent 
to $n(n+2m+3)=2^{u+1}$, which is impossible in view of the fact that    
one among numbers  $n$ and  $n+2m+3$ is an odd integer.

If $a=(2^{u+1}-1)2^u$ for a  positive integer $u$, then the equality
$S_{n,m}(2,3)=a$ is equivalent to 
 $$
n(n+2m+3)=(2^{u+1}-1)2^{u+1}.
  $$
If $2^{u+1}-1$ is a composite number, then it can be written as 
a product $2^{u+1}-1=tv$ with odd integers $t\ge 3$ and $v\ge 3$. Then 
the above equality holds for $n=v\ge 3$ and 
$m=(t(tv+1))-v-3)/2=((t^2-1)v+t-3)/2\ge 12$. 
If $2^{u+1}-1$ is a prime number, then  easily follows 
that the above equality holds only  for $n=1$ and 
$m=(2^{u+1}-1)2^{u}-2$. 

Now consider the last  case, i.e.,  when 
$a=(2^{u+1}+1)2^u$ for a  positive integer $u$. Then the equality
$S_{n,m}(2,3)=a$ is equivalent to 
 $$
n(n+2m+3)=(2^{u+1}+1)2^{u+1}.
  $$
If $2^{u+1}+1$ is a composite number, then it can be written as 
a product $2^{u+1}+1=tv$ with odd integers $t\ge 3$ and $v\ge 3$. Then 
the above equality holds for $n=v\ge 3$ and 
$m=(t(tv-1))-v-3)/2=((t^2-1)v-t-3)/2\ge 9$. 
If $2^{u+1}+1$ is a prime number, then 
 easily follows 
that the above equality holds only  for $n=1$ and 
$m=(2^{u+1}+1)2^{u}-2$.   

In view ot the above considerations,
we have shown that every integer $a\ge 4$ is equal to 
$S_{n,m}(2,3)$ for some integers $n\ge 2$ and $m\ge 1$
 in all the cases excluding the following ones:

1) $a$ is not a power of 2;

2) $a$ is not of the form $(2^{u+1}-1)2^u$, where
$2^{u+1}-1$ is a prime number and  

 3) $a$ is not of the form $(2^{u+1}+1)2^u$, where
$2^{u+1}+1$ is a prime number.
\end{example}

\begin{remark}
Notice that if $a=2^u(2^{u+1}+1)$ for  an integer 
$u\ge 1$, then $a=\sum_{i=1}^{2^{u+1}}i$, while  if $a=2^u(2^{u+1}-1)$ 
for  an integer $u\ge 1$, 
then  $a=\sum_{i=1}^{2^{u+1}-1}i$. These two identities 
together with Example 1.1 imply the  well known fact 
that every  integer $a>1$ which is not a power of $2$,
is a sum of two or more consecutive integers
(see, e.g., Dickson's History \cite[1, Ch. III, p. 139]{di}, 
where this result was attributed to Sylvester).  
  \end{remark}

  \begin{remark}
Note that it is well known (see, e.g., \cite[Subsections 2.2 and 2.3]{me1}) that in order 
to the so-called a {\it Mersenne number} 
$M_{u+1}:=2^{u+1}-1$ to be prime, $u+1$ must  itself be prime. 
A Mersenne number which is prime is called {\it Mersenne 
prime} (this is Sloane's sequence A000668  in \cite{sl}
corresponding to indices given by Sloane's sequence A000043).
Moreover, it is easy to show  that   
 in order to $2^{u+1}+1$ to be prime, $u+1$ must  be a power of 2.  
Such  numbers are in fact  {\it Fermat numbers} 
$F_{s}:=2^{2^s}+1$ ($s=0,1,2,\ldots$; 
this is Sloane's sequence A000215  in \cite{sl}). Fermat conjectured 
in 1650 that every Fermat number is prime and Eisenstein proposed as a 
problem in 1844 the proof that there are an infinite number 
of {\it Fermat primes} (i.e., Fermat numbers which are primes) 
(see \cite[p. 88]{r1}). However, the only known Fermat primes are $F_0=3$, 
$F_1=5$, $F_2=17$, $F_3=257$ and $F_4=65537$ 
(Sloane's sequence A019434 in \cite{sl}).    
 For more information on classical and alternative approaches to the  
Mersenne and Fermat  numbers, see \cite{jr}. 

Note that the conclusion at the end of Example 1.1 immediately yields the 
following interesting  assertion. 
  \end{remark}

\begin{proposition}
 The following two statements are equivalent:

\,$(i)$\, There are infinitely many Fermat primes or there are  infinitely many 
Mersenne primes$;$

$(ii)$ The set 
$\{S_{n,m}(2,3):\, n=2,3,\ldots;m=1,2,\ldots\}$  omits infinitely many 
positive integer values which are not powers of 2.
\end{proposition}

 \begin{example}
For the progression $A(3,5)=(2k+1)_{k=1}^{\infty}$ we have 
$S_{n,m}(3,5)=2(2m+4)$. From this it can be  easily seen  
that a positive integer $a\ge 8$ is equal  to some sum 
$S_{n,m}(3,5)$ with $n\ge 2$ if and only if $a$ is divisible by $4$
or $a$ is an odd composite integer greater than $14$ which is not 
a square of a prime.

More generally, if $q=p+2$, then $S_{n,m}(p,q)=n(n+2m+p-1)$.  
From this it follows that a positive integer $a$ is equal  to 
some sum $S_{n,m}(3,5)$ with $n\ge 2$ if and only if $a=4s$ with $s\ge (p+1)/2$
or  $a$ is an odd composite integer which can be expressed  as a product 
$n=ab$ with odd integers $a$ and $b$ such that $a\ge 3$ and $b\ge a+p-1$.
  \end{example}

From Examples 1.1 and 1.5 it follows that every integer greater than 10
can be expressed as a sum of two or more consecutive terms 
of the progression $A(2,3)$ or $A(3,5)$. Accordingly, it can be of interest
to consider a problem of representation of a positive integer as a sum
of two or more first consecutive integers in some progression 
$A(p,q)$. 
Notice that 
$$
S_2(p,q)=p+q,
 $$
and even Goldbach's  conjecture states that every 
even positive  integer greater than $2$ can be expressed as a sum of two 
primes. This famous conjecture  was proposed 
 on 7 June  1742 by the German mathematician Christian Goldbach 
in a  letter to Leonhard Euler \cite{go} (cf. \cite{di}).
 This conjecture has been shown to hold for all integers less 
than $4\times 10^{18}$, but remains unproven despite considerable 
effort.
 
 In view of the above equality, this conjecture is equivalent with 
 the following set equality:
  $$ 
\{S_2(p,q):\, p \,\, {\rm and}\,\, q \,\, {\rm are \,\, odd\,\, 
primes}\}=\{2n:\, n\in\Bbb N\setminus\{1,2\}\}.
  $$  
This fact suggests the investigations of the values of
$S_n(p,q)$ given by (1). Namely, for each positive integer $n$,
we will consider the  values 
 $$
S_{2n}(p,q)=n((2n-1)q-(2n-3)p),\leqno(4) 
  $$
where $p$ and $q$ are odd primes.
 
Using some heuristic arguments and computational results, 
 we propose the following ``weak even Goldbach  conjecture''.

 \begin{conjecture}[``weak even Goldbach  conjecture'']
For each even positive integer $a$ greater than  $2$
there exist a positive integer $n$ and odd primes $p$ and $q$ such that 
$a=S_{2n}(p,q)$$;$ or equivalently, that 
    $$
a=n((2n-1)q-(2n-3)p).\leqno(5)
   $$
\end{conjecture}

Clearly, the following conjecture is  stronger than Conjecture 1.6.

 \begin{conjecture}
For any  positive integer $n>1$ there exist odd primes $p$ and $q$ such that 
  $$
(2n-1)q-(2n-3)p=2.\leqno(6)
  $$
 \end{conjecture}

Note that the equality  (6) can be written as  
  $$
q=p-\frac{2(p-1)}{2n-1},
  $$
whence it follows  that  $p=2k(2n-1)+1$ and $q=2k(2n-3)+1$
for a positive integer $k$. Hence, Conjecture 1.7 is equivalent to the 
following one.

\vspace{2mm}

\noindent{\bf Conjecture 1.7'.} {\it For  any  integer $n>1$ there exists 
a positive integer $k$ such that both numbers 
$p=2k(2n-1)+1$ and $q=2k(2n-3)+1$ are primes.}
\vspace{2mm}

If $p$ and $q$ are odd primes, then from the expression (1) we see that 
$S_n(p,q)$ is odd if and only if $n$ is even. The following conjecture
is the odd analogue of Conjecture 1.6.

 \begin{conjecture}[``weak odd Goldbach  conjecture'']
For each odd positive integer $a$ greater than  $2$
there exist a positive integer $n$ and odd primes $p$ and $q$ such that 
$a=S_{2n+1}(p,q)$$;$ or equivalently, that 
    $$
a=(2n+1)(nq-(n-1)p).\leqno(7)
   $$
 \end{conjecture}

Clearly, the following conjecture is  stronger than Conjecture 1.8.

\begin{conjecture}
For any  positive integer $n>1$ there exist odd primes $p$ and $q$ such that
  $$
nq-(n-1)p=1.\leqno(8)
  $$
 \end{conjecture}

 From the equality   (8) we have 
  $$
q=p-\frac{p-1}{n},
  $$
whence we conclude that  $p=nk+1$ and $q=(n-1)k+1$ for a positive 
integer $k$. This together with the fact that $k=p-q$ is even shows that 
Conjecture 1.9 is equivalent to the following one.

\vspace{2mm}

\noindent{\bf Conjecture 1.9'.} {\it For  any  integer $n>1$ there exists 
a positive integer $k$ such that both numbers 
$p=2kn+1$ and $q=2k(n-1)+1$ are primes.}
\vspace{2mm}

Finally, notice that Conjectures 1.6 and 1.8 can be 
joined into the following conjecture.

 \begin{conjecture}[``weak Goldbach  conjecture'']
Conjectures $1.6$ and $1.8$ are  true if and only if 
the following statement holds true$:$

For each  positive integer $a$ greater than  $2$
there exist a positive integer $n$ and odd primes $p$ and $q$ such that 
    $$
a=\frac{n}{2}((n-1)q-(n-3)p).\leqno(9)
   $$
\end{conjecture}

\end{document}